\documentclass[review]{elsarticle}

\usepackage{lineno,hyperref}
\usepackage{graphicx}
\usepackage{amsmath}
\usepackage{amssymb}
\usepackage{float}
\usepackage{bbm}

\modulolinenumbers[5]

\def\CoVaR{\mathop{\operatorname{CoVaR_{\alpha,\beta}^{=}}}}

\def\VaRa{\mathop{\operatorname{VaR_\alpha}}}
\def\VaRb{\mathop{\operatorname{VaR_\beta}}}
\def\quantA{a}
\def\quantAA{a\,}
\def\quantB{b}
\def\quantBB{b\,}
\def\Newopt{f}
\def\funH{g}
\def\rootfun{s}
\def\solFun{h}
\def\nowa{F}

\newtheorem{lemma}{Lemma}
\newtheorem{theorem}{Theorem}

\journal{Journal of \LaTeX\ Templates}

\biboptions{authoryear}

\begin{document}

\begin{frontmatter}

\title{CoVaR-based portfolio selection 
 }


\author[mymainaddress,mysecondaryaddress]{Anna Patrycja Zalewska \corref{mycorrespondingauthor}}

\cortext[mycorrespondingauthor]{Correspondence to:}
\ead{A.Zalewska@mini.pw.edu.pl}
\ead{A.Zalewska@mimuw.edu.pl}

\address[mymainaddress]{Faculty of Mathematics, Informatics and Mechanics, \\ University of Warsaw, Warsaw 02097, Poland }
\address[mysecondaryaddress]{Faculty of Mathematics and Information Science, 
 \\Warsaw University of Technology, Warsaw 00662, Poland}

\begin{abstract}
We consider the portfolio optimization with risk measured by conditional   \\ \mbox{value-at-risk}, based on the stress event of chosen asset being equal to the opposite of its \mbox{value-at-risk} level, under the normality assumption. Solvability conditions are given and illustrated by examples. 
  \end{abstract}


\end{frontmatter}

\linenumbers

\section{Introduction}

What is the best asset allocation? 
\cite{M52} answered this question in his~landmark article on~\mbox{mean-variance} model of~portfolio selection, followed few years later by~a~book---
\cite{M59}---on the same subject and thus the Modern Portfolio Theory originated, providing its creator with Alfred Nobel Memorial Prize in Economic Sciences~\mbox{(1990)}\footnote{Awarded jointly to Harry M. Markowitz, Merton H. Miller and William F. Sharpe. }. Markowitz exceptional idea was to~minimize portfolio variance or~(equivalently) standard deviation for a~fixed portfolio return $E$. 
Therefore, the investor faces an~optimization problem, given as~follows:
\begin{equation}\label{StDevOpt} \left\lbrace \begin{array}{lcccccc}
			\operatorname{\sigma}(x) \rightarrow \min \\ 
			 x^{\mathrm{T}} \mu \geq E\\
			 x_1 + \dots + x_n = 1 \\
			 x_1, x_2, \dots, x_n \geq 0 \quad \left(\star \right)\\
         \end{array} \right. 
          \end{equation}
where portfolio (i.e.~investment strategy) \mbox{${x = (x_1, \dots, x_n)^\mathrm{T} }$} meets the natural condition of~summing~up~to~$1$ 
and~${\mu = (\mu_1, \dots, \mu_n)^\mathrm{T}}$ is defined as expected value of~\mbox{$n$-dimensional} random variable of returns on \textbf{risky} assets, ${R=(R_1, \dots, R_n)^\mathrm{T}}$, with~${X = x^{\mathrm{T}} R=\sum\limits_{i = 1}^n x_iR_i}$ being the univariate random variable of~return on the~portfolio. The~$\left(\star \right)$ constraint is optional and~concerns the~possibility of \mbox{short-selling}. Without it the~model is~known as~\mbox{the Black model}---see~\cite{PA} and the~original article of~\cite{Black72}. 

 To obtain a~non-degenerate problem, two assumptions are made. First, ${\mu \nparallel \mathbbm{1}_n =(1,\dots, 1)^\mathrm{T}}$, i.e.~not every asset has the~same expected return. Second, the~covariance matrix of~$R$, ${\Sigma = \left[ \sigma_{ij}^{} \right]_{i,j = 1,\dots, n} }$ is positive definite.
 
Merton solved problem~\eqref{StDevOpt} with those two assumptions satisfied and constraint~$\left( \star \right)$ dropped. 
His formulae in contemporary terms~(\mbox{cf.~\cite{merton}}) assume the~following form:
\begin{equation}\label{merton}
\operatorname{x}(E) = \Big(\big(\alpha_{M}^{} \gamma_{M}^{} - \beta_{M}^2\big)\Sigma\Big)^{-1} \left(\, 
 \left| \begin{array}{ccccccc}
			E & \beta_{M}^{} \\
			 1 & \gamma_{M}^{}
          \end{array}\right| \mu + 
           \left| \begin{array}{ccccccc}
			\alpha_{M}^{} & E \\
			\beta_{M}^{} & 1
          \end{array}\right| \mathbbm{1}_n
 \right)
\end{equation}
where
$ \alpha_{M}^{} = \mu^\mathrm{T} \Sigma^{-1} \mu,\ \beta_{M}^{} = \mu^\mathrm{T} \Sigma^{-1}\mathbbm{1}_n^{},\ \gamma_{M}^{} = \mathbbm{1}_n^\mathrm{T} \Sigma^{-1} \mathbbm{1}_n^{} $. 
The graph of ${x(E)}$ is a~line~$\Gamma$ called the~\mbox{\textbf{critical line}}. Portfolio is said to~be~\textbf{efficient} if there is no portfolio with either smaller $\sigma$ for the same or greater $E$, or greater $E$ for the same or smaller $\sigma$. Both terms first appear in~\cite{M52}. The set of~efficient portfolios is a~subset of~${\lbrace x(E) \mid E\in \mathbb{R} \rbrace}$. Its image in~${x \mapsto (\sigma(x),E(x))}$ mapping is~known as the~\textbf{efficient frontier}.

The main objective of~this paper is to~provide some insight into problem~\eqref{StDevOpt} with $\CoVaR$ as the~alternate risk measure, without $\left( \star \right)$ constraint (if~not stated differently), and under normality assumption added to original ones. We begin with briefly stating the reasons why $\VaRa$ is not of~interest in~that case and then proceed with $\CoVaR$-a conditional value at risk proposed by~\cite{adrian2008covar,adrian2016covar}, not~to~be~confused with $\operatorname{CVaR}$ (Mean Excess Loss) as used by~\cite{rockafellar2000optimization} for the~optimization problem. With the aid of~examples the properties of~the~new risk measure, the critical set i.e. the set of {minimum-$CoVaR$} portfolios for fixed expected value and the very existence of that set are discussed.

\section{Risk measured by VaR}
For $X$-a~random variable of~an~asset portfolio (the \mbox{profit/loss} approach) let value-at-risk be defined as ${\VaRa (X) = -Q^{+}_{\alpha}(X)}$.
It is worth noting that $\operatorname{VaR}$ is a downside risk measure, while $ \operatorname{\sigma}$ is classified as a volatility measure. The former is monotone, translation invariant and positively homogeneous, lacking only subadditivity to be a coherent risk measure (\mbox{cf.~\cite{artzner1999coherent}}), while the latter is just positively homogeneous. However, for ${\alpha \in ( 0, 1\slash 2)}$ and under normality assumption, $\VaRa$ is~coherent, as~\cite{artzner1999coherent} prove. It might seem promising, but as soon as we calculate the actual \mbox{value-at-risk} of portfolio, 
${-x^{\mathrm{T}} \mu + \operatorname{\sigma}(x) \cdot (- \Phi^{-1}(\alpha))}$, we can clearly see that it~yields the same solutions as~$\sigma(x)$ (for formal proofs see \cite{alexander2002economic}).

\section{Risk measured by $\mathrm{CoVaR}^{=}$}\label{sec:CoVaR} 
In present section we begin by giving a definition of $\CoVaR$, as introduced by~\cite{adrian2008covar,adrian2016covar}, though notation is rather that of~\cite{MaiScha} (cf.~\cite{bernardi2016covar}).
\begin{equation*}
\CoVaR(X \mid Y) = \VaRb(X \mid Y = -\VaRa (Y))
\end{equation*} 

The first assumption to be made is that of normality, $ {R \sim \mathcal{N}(\mu, \Sigma)}$, ${\Sigma >0}$. 
For bivariate Gauss distribution, where
\begin{equation*} 
 \left( \begin{array}{rrr}
 X  \\
 Y \\
 \end{array}\right)
 \sim 
 \mathcal{N}_2 \left( 
\left( \begin{array}{ccccccc}
\mu_X^{} \\ 
\mu_Y^{}\\
  \end{array}\right) 
 ,
   \left( \begin{array}{cccc}
 \sigma_X^2 & \rho_{X,Y}^{}\, \sigma_X^{} \sigma_Y^{}\\
 \rho_{X,Y}^{}\, \sigma_X^{} \sigma_Y^{} &  \sigma_Y^2 \\
 \end{array}\right)
 \right) \quad \mbox{for } |\rho_{X,Y}^{}|<1
 \end{equation*}
and ${X = \mu_X^{}  \pm \sigma_X^{} \slash \sigma_Y^{} (Y- \mu_Y^{})}$ 
for ${\rho_{X,Y}^{} = \pm 1}$, 
\cite{MaiScha} obtain formula:
\begin{equation}\label{CoVaR=Gaussformula}
\CoVaR(X\mid Y) = -\mu_X^{} - \sigma_X^{} \Big( \rho_{X,Y}^{} \Phi^{-1} (\alpha) + \Phi^{-1}(\beta) \sqrt{1-\rho_{X,Y}^2}\Big)
\end{equation}

The second assumption, ${\alpha,\beta \in (0,1 \slash 2) }$, is only natural as investor interest in calculating $\operatorname{VaR}$ lays chiefly in significance level being close to $0$. Consequently, ${a = -\Phi^{-1}(\alpha)}$ and ${b= -\Phi^{-1}(\beta)}$ are positive numbers which prevents us from dwelling~on \mbox{sub-cases}.

This established, the distribution of $X=x^{\mathrm{T}}R$ is conditioned on one chosen variable $R_i$, $i \in \lbrace1, \dots, n \rbrace$. Without loss of generality let that be $Y = R_1$. Naturally, $\VaRa(R_1) = -\mu_1 + \quantAA \sigma_1 $ and 
${ \rho_{X,Y} = 
  \sigma_1^{-1} \big( x^{\mathrm{T}} \Sigma x \big)^{-1\slash 2} \sum\limits^{n}_{k=1} x_i^{} \sigma_{1k} }$.

Investor faces the following optimization problem:

\begin{equation}\label{CoVaR=} \left\lbrace \begin{array}{lcccccc}
			\CoVaR(X \mid Y) \rightarrow \min \\ 
			 x^{\mathrm{T}} \mu = E\\
			 x_1 + \dots + x_n = 1 \\
          \end{array} \right. 
          \end{equation}
For ${ \rho_{X,Y} = \pm 1}$ there is a linear relationship between $X$ and~$Y$. 
Consequently ${X \mid Y = -\VaRa(Y)}$ is a constant, hence respectively ${\CoVaR(X \mid Y) = -\mu_1 \pm \quantAA \sigma_1} $.

We observe (after applying  formulae ${\sigma_{ij} =\rho_{ij}\sigma_i, \sigma_j}$ and \eqref{CoVaR=Gaussformula}, and a~little manipulation) that the following is true: 
\begin{equation}\label{covFormula}
\CoVaR(X \mid Y ) =  -  x^{\mathrm{T}}_{}\mu + \quantA \sum\limits^{n}_{k=1} x_i^{} \rho_{1k}^{} \sigma_k^{}  + \quantBB \sqrt{ x^{\mathrm{T}}_{} Q x }
\end{equation} 
where $Q$ is defined as:
\begin{equation*} Q =   \Sigma  - q \cdot  q^{\mathrm{T}} = \newcommand*{\temp}{\multicolumn{1}{r|}{}}
\left(\begin{array}{ccccc}
0 \mskip-10mu & \temp  & 0 \\ \cline{1-3}
0  \mskip-10mu  & \temp  & \widehat{Q}  \\
\end{array}\right)
\end{equation*}
with 
 $\widehat{Q}$ being a~symmetric $(n-1) \times (n-1) $ positive definite matrix (see the Appendix) and   ${q = 1\slash\sigma_1 \Sigma\, e_1}$
where ${e_1 = (1, 0, \dots, 0)^{\mathrm{T}}}$.

Since ${\CoVaR(X \mid Y)}$ depends solely on $x$ for given ${\Sigma, \mu, \alpha, \beta}$, let~us from now~on denote is as ${\CoVaR(x)}$. Therefore  optimization problem~\eqref{CoVaR=} presents itself as follows:
\begin{equation}\label{CoVaR=2} \left\lbrace \begin{array}{lcccccc}
\CoVaR(x)= - x^{\mathrm{T}}\mu +  \quantAA  x^{\mathrm{T}} q +\quantBB \sqrt{ x^{\mathrm{T}} Qx  }\rightarrow \min  \\ 
			  x^{\mathrm{T}}  \mu =  E\\
			 x_1 + \dots + x_n = 1  \\
          \end{array} \right. 
          \end{equation}
Obviously for ${\quantA = \quantB}$ the critical set remains independent of~$\quantA$. 
Function $\CoVaR $ is convex (see Appendix) and positive homogeneous (of~degree~$1$), as $\operatorname{\sigma}$~is. 

$\CoVaR$ is not  bounded above and it does not have to be bounded below\footnote{E.g. should partial derivative of $\CoVaR$ with respect to $x_1$, i.e. ${-\mu_1 + a \sigma_1}$ be a~\mbox{non-zero} number, the function is not bounded below. On the other hand, consider an example with diagonal $\Sigma$ and ${\mu = \mu_1 e_1}$ where ${\mu_1 = a \sigma_1}$. Then function is bounded below.} (contrary to~$\operatorname{\sigma}$ which is always bounded below).

What conditions should be met in~order for problem~\eqref{CoVaR=2} to~have a~solution? Before introducing the~main~theorem of~this work, we define  ${\widehat{x} = (x_2,\dots, x_n)^{\mathrm{T}}}$, ${\Newopt(\widehat{x}) = \CoVaR(1-x_2-\dots-x_n,x_2, \ldots, x_n)}$, ${\widehat{\mu} = (\mu_2-\mu_1,\dots, \mu_n - \mu_1)^{\mathrm{T}}}$, \\${\widehat{q} = (\rho_{12}\sigma_{2}- \sigma_1, \dots,  \rho_{1n} \sigma_{n} - \sigma_1)^{\mathrm{T}}}$, ${\widehat{E} = E-\mu_1}$ and
\begin{align*}
G = 
\left(\begin{array}{ccccc}
 \alpha_{C}^{} & \beta_{C}^{} \\
 \beta_{C}^{} & \gamma_{C}^{}  \\
\end{array}\right)
\ \mbox{where } \
\alpha_{C}^{} = \widehat{\mu}^{\mathrm{T}} \widehat{Q}^{-1} \widehat{\mu}, \ 
\beta_{C}^{} = 
 \widehat{\mu}^{\mathrm{T}} \widehat{Q}^{-1} \widehat{q}
\ \mbox{and} \
\gamma_C^{} = \widehat{q}^{\mathrm{T}} \widehat{Q}^{-1} \widehat{q} 
\end{align*}
with  ${\Delta=b^2\alpha_C^{} - a^2\det G }$.
\begin{theorem}\label{th1}
Let vectors ${ \mathbbm{1}_n, \mu, q}$ be linearly independent. 
\\If ${\Delta>0}$ then the  optimization problem 
\begin{equation}\label{CoVaR=3} \left\lbrace \begin{array}{lcccccc}
\Newopt(\widehat{x})= -\mu_1 + \quantAA \sigma_1 -\widehat{x}^{\mathrm{T}}\widehat{\mu} +  \quantAA \widehat{x}^{\mathrm{T}} \widehat{q} +\quantBB \sqrt{ \widehat{x}^{\mathrm{T}} \widehat{Q}\widehat{x}  }
\rightarrow \min  \\ 
			  \widehat{x}^{\mathrm{T}}  \widehat{\mu} =  \widehat{E}\\
          \end{array} \right. 
          \end{equation}
equivalent to~problem~\eqref{CoVaR=2} has for a~given $\widehat{E}$ a~unique solution
\begin{align*}
\widehat{x}( \widehat{E}) &= \dfrac{ \widehat{E}}{\alpha_{C}^{}}  \, \widehat{Q}^{-1}\widehat{\mu}
+
| \widehat{E}| \,  \dfrac{\quantAA }{\alpha_{C}^{} \sqrt{\Delta}} \,\widehat{Q}^{-1}(\beta_C^{}\, \widehat{\mu} - \alpha_C^{}\,
     \widehat{q}) 
      ,
\end{align*}
and 
\begin{align*}
\Newopt\Big(\widehat{x}( \widehat{E})\Big) =  
   -\mu_1 + \quantAA \sigma_1 + \widehat{E} \left( \frac{\quantAA \beta_C^{} }{ \alpha_C^{}}- 1\right) + \frac{|\widehat{E}| }{ \alpha_C^{}} \sqrt{ \Delta }.
 \end{align*}
 Moreover, for ${\Delta>0}$ and ${x(\widehat{E}) = (1- \mathbbm{1}_{n-1}^{\mathrm{T}} \widehat{x}(\widehat{E}), \widehat{x}(\widehat{E})^\mathrm{T})^\mathrm{T}} $ the~following is~true:
 \begin{enumerate}
 \item For ${ \quantA \beta_C^{}- \alpha_C^{} \leq - \sqrt{\Delta} }$ there is no \mbox{$\CoVaR$-efficient portfolio.}
 \item For ${ - \sqrt{\Delta}<\quantA \beta_C^{}- \alpha_C^{} \leq \sqrt{\Delta} }$  only ${x( \widehat{E})}$ for~${\widehat{E}\geq 0}$ are \mbox{$\CoVaR$-efficient} portfolios.
 \item For ${ \sqrt{\Delta}<  \quantA \beta_C^{} -\alpha_C^{} }$ all portfolios ${x( \widehat{E})} $ constitute the set \\of~\mbox{$\CoVaR$-efficient} portfolios.
 \end{enumerate}
\end{theorem}
 We relegate the proof to Appendix.

Now some examples will be presented. The first shows how unless ${\Delta > 0}$ condition is satisfied  we obtain optimization problem without solution
and suggests using ${(\star)}$ constraint in such a~case. 

\subsection{Example 1}\label{MinNotAttained} 

\vspace{-\baselineskip}
\begin{equation*}{(R_1, R_2,R_3) \sim \mathcal{N}_3 \left( 
\left( \begin{array}{ccccccc}
1 \\ 
4 \\
 3 \\
  \end{array}\right) 
 ,
   \left( \begin{array}{rrr}
1 &  -4\slash 3 &  2\slash 3 \\ 
 -4\slash 3 &  4 &  -1 \\
 2\slash 3 &  -1 &  1 \\
 \end{array}\right)
 \right),
 \quad \quantA = 8\slash 10, \ \quantB =  7\slash 10}
  \end{equation*}
By a straightforward calculation we get: 
\begin{equation*}
{\CoVaR(x) = 1\slash 30 \Big(-6 x_1^{} -152 x_2^{} - 74 x_3^{} + 7 \sqrt{20 x_2^2 - 2 x_2^{} x_3^{} + 5 x_3^2} \,\Big)} 
\end{equation*}
and  ${\hat{Q} =  \left( \begin{array}{rrr}
 20\slash 9 & -1\slash 9\\
 -1\slash 9 &  5\slash 9\\ \end{array}\right)}.$
Suppose ${x^{\mathrm{T}} \mu = 2}$---together with ${x_1 + x_2 + x_3 =1}$ condition it implies ${x_2 = 2 x_1-1}$ and ${x_3 = - 3 x_1+ 2}$. Without~$\left( \star \right)$ constraint we can solve the optimization problem~\eqref{CoVaR=2} for ${E = 2}$ by~minimizing over $\mathbb{R}$ the following function:
\begin{align*}
\funH (x_1^{}) &= \CoVaR(x_1^{}, 2 x_1^{}-1, - 3 x_1^{}+ 2)= \\
&=  1\slash 30 \Big(-88 x_1^{}+4 +7 \sqrt{137 x_1^2-154 x_1^{}+44} \,\Big)
\end{align*}
Clearly, ${\lim\limits_{x_1 \rightarrow + \infty}\funH(x_1) = - \infty}$, hence $ \funH$ does not attain minimum, which shows that even with additional constraints in this simple example $\CoVaR$ may be unbounded below. 
Now we add the \mbox{non-negativity} constraint~$(\star)$ i.e. we minimize~$\CoVaR$ on the standard \mbox{$2$-simplex} (equilateral triangle in~$\mathbb{R}^3$): 
\begin{equation*} \left\lbrace \begin{array}{lcccccc}
g(x_1) \rightarrow \min  \\ 
			 x_1 \in [1\slash 2, 2\slash 3]\\
          \end{array} \right. \quad  \mbox{ gives us }x_{\min} = 2\slash 3,\ g(x_{\min}) = 1\slash 45 (-82 + 7 \sqrt{5}).
          \end{equation*}
          
\begin{figure}[!ht]
\centering
\includegraphics{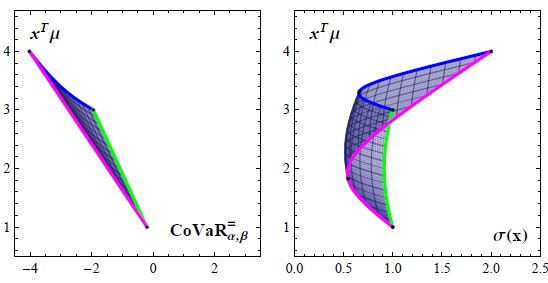} 
\vspace{-1.5\baselineskip}
\caption{ Images of functions  ${x \mapsto \big(\CoVaR(x), x^\mathrm{T} \mu \big)}$ (left) and ${x \mapsto \big(\operatorname{\sigma}(x), x^{\mathrm{T}} \mu \big)}$ (right) for ${x \in \lbrace (x_1, x_2,x_3) \in \mathbb{R}^3 \mid  x_1+x_2+x_3 =1,\ x_1 \geq 0,\ x_2 \geq 0,\ x_3 \geq 0 \rbrace }$. } 
\label{ex1_01}
\end{figure}
Figure~\eqref{ex1_01} might be of~use both as an illustration and providing comparisons.  

\bigskip
With one example without the counterpart of the critical line~\eqref{merton} the second one will be presented---this time with function $\CoVaR$ attaining its minimum over $ \lbrace x \in \mathbb{R}^n\mid  x_1 + \dots + x_n =1\rbrace$ which allows us to~consider~\eqref{CoVaR=2} both with and without the~non-negativity constraint~$(\star)$. 

\subsection{Example 2}\label{MinAttained}
\vspace{-\baselineskip}
\begin{equation*} (R_1, R_2,R_3) \sim \mathcal{N}_3 \left( 
\left( \begin{array}{ccccccc}
2 \\ 
3 \\
 1 \\
  \end{array}\right) 
 ,
   \left( \begin{array}{rrr}
 1 &  1\slash 5 &  1 \\
 1\slash 5 &  1 &  0 \\
   1 &  0 &  9 \\
 \end{array}\right)
 \right), \quad \quantA =1, \ \quantB= 2
 \end{equation*}
 
In this case ${\CoVaR(x) = 1 \slash 5 \Big(-5 x_1^{} - 14 x_2^{} + 2 \sqrt{24 x_2^2 - 10 x_2^{} x_3^{} + 200 x_3^2 \ }\Big)}$. 
Using just the \lq portfolio constraint\rq, $x_1 + x_2 + x_3=1$ from the optimization problem~\eqref{CoVaR=2} generates a bounded below function with unique global minimum:
\begin{equation*} \left\lbrace \begin{array}{lcccccc}
\CoVaR(x) \rightarrow \min  \\ 
			 x_1 + x_2 + x_3 = 1\\
          \end{array} \right. \quad  \mbox{ gives us }x_{\min} = (1,0,0)^{\mathrm{T}},\ \CoVaR(x_{\min}) = -1.
          \end{equation*}
\subsection{Example 3}
\begin{figure}[!ht]
\centering
\includegraphics[width=.3\textwidth]{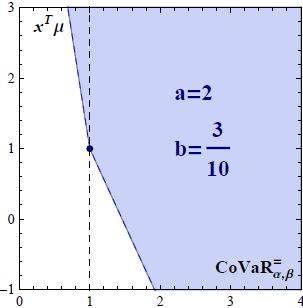}\quad
\includegraphics[width=.3\textwidth]{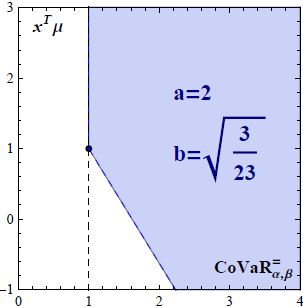}

\medskip

\includegraphics[width=.3\textwidth]{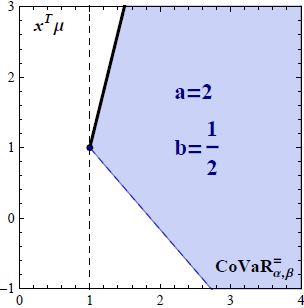} \quad 
\includegraphics[width=.3\textwidth]{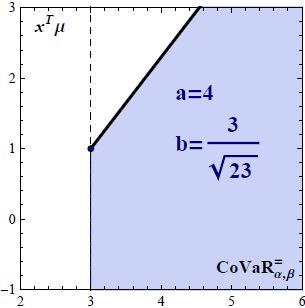}\quad
\includegraphics[width=.3\textwidth]{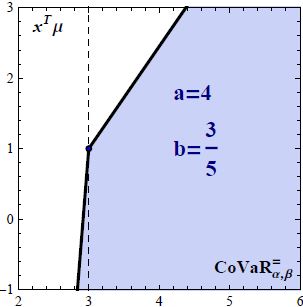}

\caption{Images of ${x \mapsto \big(\CoVaR(x), x^\mathrm{T} \mu \big)}$ for a~choice of~$(\quantA, \quantB)$.  \newline 
The~$\CoVaR$-efficient frontiers are marked by thick lines. }
\label{ex4_01}
\end{figure}

\vspace{-\baselineskip}
\begin{equation*}  (R_1, R_2,R_3) \sim \mathcal{N}_3 \left( 
\left( \begin{array}{ccccccc}
1 \\ 
2 \\
3 \\
 \end{array}\right) 
 ,
   \left( \begin{array}{ccccccc}
1 & 1 & 2\\
1 & 9 & 0\\
2 & 0 & 16\\
 \end{array}\right)
 \right)
\end{equation*}      
where:
\begin{equation*} \widehat{Q} = 
\left( \begin{array}{ccccccc}
8 & -2\\
-2 & 12\\
 \end{array}\right) 
\end{equation*}  
This example is essentially an~illustration of~the~dependence between values of~$\quantA, \quantB$ and the~obtained set of~$\CoVaR$-efficient portfolios. We have 
\begin{equation*}
\CoVaR(x)=(\quantA-1)x_1^{}+(\quantA-2)x_2^{}+(2 \quantA-3 )x_3^{}+2 \quantB \sqrt{2 x_2^{2}-x_2^{} x_3^{}+3 x_3^{2}\,}.
\end{equation*}
Figure~\eqref{ex4_01} presents five different possibilities, of~which two first represent case~$1$ in the theorem, two next---case~$2$ and the last one case~$3.$

\section{Appendix: proofs $\&$ auxiliary results}\label{app}

\subsection{Positive definiteness of $\hat{Q}$ and positive semidefiniteness of $Q$}
For $\Sigma>0$ and ${e_1 = \big(1,0, \dots,0 \big)^\mathrm{T} }$ there is:
\begin{equation*}
Q=\Sigma - \sigma_1^{-2} \, \Sigma e_1^{} e_1^{\mathrm{T}} \Sigma
=\Sigma \big(\Sigma^{-1}_{} - \sigma_1^{-2} \,  e_1^{} e_1^{\mathrm{T}} \big) \Sigma.
\end{equation*}
Matrix ${e_1^{} e_1^{\mathrm{T}}}$ has only one non-zero element (i.e. $1$ in first row, first column). In~consequence ${\Sigma^{-1}_{} - \sigma_1^{-2} \,  e_1^{} e_1^{\mathrm{T}} }$ has rank equal at least ${n-1}$ as its columns from second to last are the corresponding columns of~$\Sigma^{-1}_{}$. Therefore matrices $\hat{Q}$ and $Q$ are of the rank  ${n-1}$, so that they are positive definite and positive semidefinite, respectively.
\\
On the other hand, for $X=x^{\mathrm{T}}R$ and ${Y = R_1}$ we have
\begin{multline*}
\forall  x \in \lbrace (x_1, \widehat{x}) \in \mathbb{R}^n \mid  x_1 \in \mathbb{R},\ \widehat{x} \in \mathbb{R}^{n-1} \rbrace : \\
 \widehat{x}^{\mathrm{T}\, } \widehat{Q} \widehat{x}  = x^{\mathrm{T}\, }  Q x = x^{\mathrm{T}\, }  \Sigma x -  x^{\mathrm{T}\, } q\, q^{\mathrm{T}\,} x = \sigma^2_{X} (1 - \rho_{X,Y}^2) \geq 0
\end{multline*}
so that ${\rho_{X,Y}^2 = \pm 1}$ if and only if ${x = (const, 0, \dots, 0)^{\mathrm{T}}}$ (due to positive definiteness of $\widehat{Q} $) which in problem~\eqref{CoVaR=} implies ${x = (1, 0, \dots, 0)^{\mathrm{T}}}$.

\subsection{Convexity of $\CoVaR$ function}

Function $\rootfun(\widehat{x}) = \sqrt{\widehat{x}^{\mathrm{T}\, } \widehat{Q} \widehat{x}}$
\, is a~norm induced by~inner product ${\widehat{x}^{\mathrm{T}} \widehat{Q} \widehat{x} }$ and as such is convex on~$\mathbb{R}^{n-1}$--indeed,
due to~$\rootfun$ being positively homogeneous of~degree~$1$, the~Jensen's inequality can be rewritten
as:
\begin{equation}\label{triangle}
\forall \widehat{x} ,\widehat{y} \in \mathbb{R}^{n-1}, \forall \lambda \in [0, 1] :
\quad
 \rootfun \big(\lambda \widehat{x}+(1-\lambda)\widehat{y} \big) \leq  \rootfun(\lambda\widehat{x})+\rootfun\big((1-\lambda) \widehat{y}\big)
\end{equation}  
which is just triangle inequality, satisfied by~every norm.  

The fact that ${\rootfun(x_2, 0,\dots,0) = \vert x_2 \vert \sqrt{\widehat{q}_{11}}}$ contradicts the~strict convexity. However, as the inequality~\eqref{triangle} is strict provided that ${\widehat{x}   \nparallel \widehat{y}}$, function  $\rootfun$ is strictly convex on~every line segment not contained in any \mbox{half-line} with the origin as its element.

The instantaneous implication is that 
${x^{\mathrm{T}} Q x = \widehat{x}^{\mathrm{T}} \widehat{Q} \widehat{x}}\,$ is convex on $\mathbb{R}^n$. 
Function
${\CoVaR(x)= - x^{\mathrm{T}} \mu + \quantAA x^{\mathrm{T}} q  + \quantBB \sqrt{ x^{\mathrm{T}} Qx  }}$ is convex as a~linear combination of~(proper) convex functions with positive coefficients.

\subsection{A lemma about certain convex function}\label{lemma}
\begin{lemma}
Let us consider function ${\nowa:\mathbb{R}\mapsto \mathbb{R}}$, ${\nowa(t) = s\cdot t + \sqrt{(t - p)^2 + q\,}}$ with ${s \in \mathbb{R}_+ \cup \lbrace 0 \rbrace,q \in \mathbb{R}_+,\ p \in \mathbb{R}}$. Then ${\nowa}$ is convex and the following is true:
 \begin{enumerate}
\item For ${s \in [0,1)}$ there is a global minimum ${p s + \sqrt{q (1 - s^2)\,}}$ attained\\ at~${t = p - s \sqrt{q(1 - s^2)^{-1}\,}}$
\item For ${s=1}$ there is no global minimum, but ${\nowa(t)>\lim\limits_{\xi \rightarrow -\infty}\nowa(\xi)=p }$
\item For $s>1$ function $\nowa$ is not bounded below, ${\lim\limits_{\xi \rightarrow -\infty}\nowa(\xi)= - \infty }$.
\end{enumerate} 
\end{lemma}
Proof via direct calculation.

\subsection{Proof of the main theorem}
For ${\mu, \mathbbm{1}_n, q}$ linearly independent (which implies ${n \geq 3}$) we solve the following problem:
\begin{equation}\label{CoVaR=4}
 \left\lbrace \begin{array}{lcccccc}
\widehat{\solFun}(\widehat{x})= 1/2\, \widehat{x}^{\mathrm{T}} \widehat{Q} \widehat{x} \rightarrow \min  \\ 
			  \widehat{x}^{\mathrm{T}} \widehat{\mu} = \widehat{E}\\
			  \widehat{x}^{\mathrm{T}} \widehat{q} = \widehat{t}\\
\end{array} \right. 
\end{equation}
and then vary the~parameter~$\widehat{t}$ in order to~minimize the~obtained solution with respect to~$\widehat{t}$. Naturally, $\widehat{\solFun}$ is strictly convex. 

Via~the~method  of~Lagrange multipliers we get ${\widehat{Q}\widehat{x} = \lambda_1 \widehat{\mu} +\lambda_2 \widehat{q} }$. In~consequence ${\widehat{x} = \lambda_1 \widehat{Q}^{-1}\widehat{\mu}  +\lambda_2 \widehat{Q}^{-1} \widehat{q} }$. Left multiplication of first equation (by $\widehat{x}^{\mathrm{T}}$) and the second equation (by $ \widehat{\mu}^{\mathrm{T}}$ and $\widehat{q}^{\mathrm{T}}$)  gives us:
\begin{align}
\label{eq1}
\widehat{x}^{\mathrm{T}} \widehat{Q} \widehat{x} = & \lambda_1 \widehat{E} + \lambda_2 \widehat{t}
\\ \label{eq2}
\widehat{E}= & \lambda_1 \,\alpha_C^{}  +\lambda_2\,  \beta_C^{}
\\ \label{eq3}
\widehat{t} = & \lambda_1 \, \beta_C^{}  +\lambda_2\, \gamma_C^{}
\end{align}
Scalars ${\lambda_1, \lambda_2}$ are easily obtained from \eqref{eq2} and~\eqref{eq3}. 
Ultimately, 
\begin{equation*}
{\widehat{x}^{\mathrm{T}} \widehat{Q} \widehat{x}  = \Big(\widehat{E}, \widehat{t} \, \Big) \cdot  G^{-1} \cdot \Big(\widehat{E}, \widehat{t} \, \Big)^{\mathrm{T}} },
\end{equation*}where ${G}$ 
is the Gramian matrix of linearly independent vectors  $\widehat{\mu}, \widehat{q} $, with the inner product defined by matrix $\widehat{Q}^{-1} $ ($\widehat{Q}$ being positive definite implies positive definiteness of both $G$ and $G^{-1}$).
Unless ${\widehat{E}=0}$ this quadratic function of~$\widehat{t}$ is positive and bounded  from $0$ (in the former case $0$ is attained for ${\widehat{t}=0}$ and equation \eqref{eq1} yields ${\widehat{x}=0}$). Then,
\begin{equation*}
\Newopt(\widehat{x})= -\mu_1 + \quantAA \sigma_1 - \widehat{E} +  \quantAA   \widehat{t} +\quantBB \big(\alpha_C^{}\slash \det G \big)^{1\slash 2} \ \sqrt{ \left( \widehat{t} - \beta_C^{}\slash \alpha_C^{ } \, \widehat{E} \, \right)^2  +  \det G \cdot  \widehat{E}^2 \slash  \alpha_C^{2} \, }
\end{equation*}

We take
${ \quantB^{-1} \big(\alpha_C^{}\slash \det G \big)^{-1\slash 2}  \cdot \Big( \Newopt(\widehat{x})  +\mu_1 - \quantAA \sigma_1 + \widehat{E} \Big) }$. Then the~obtained function is of~the type from~lemma~\ref{lemma} with ${s = \quantA \quantB^{-1} \big(\det G \slash \alpha_C^{} \big)^{1\slash 2} }$. That means the~global minimum is achieved for 
\begin{equation*} \widehat{t }
= \dfrac{ \beta_C^{}}{\alpha_{C}^{}} \widehat{E}- 
| \widehat{E}| \,  \dfrac{\quantAA \det G }{\alpha_{C}^{} \sqrt{\Delta}}  ,
\end{equation*}
 under condition~${b > a \big(\det G \slash \alpha_C^{} \big)^{1\slash 2}}$ (\mbox{i.e} ${s<1}$) equivalent to ${\Delta >0}$ as ${a,b>0}$.

 Then, ${\widehat{x} = \lambda_1 \widehat{Q}^{-1}\widehat{\mu}  +\lambda_2 \widehat{Q}^{-1} \widehat{q} }$ yields formula for ${\widehat{x}(E)}$:
\begin{align*}
     \widehat{x}( \widehat{E}) &= \dfrac{ \widehat{E}}{\alpha_{C}^{}}  \, \widehat{Q}^{-1}\widehat{\mu}
+
| \widehat{E}| \,  \dfrac{\quantAA }{\alpha_{C}^{} \sqrt{\Delta}} \,\widehat{Q}^{-1}(\beta_C^{}\, \widehat{\mu} - \alpha_C^{}\,
     \widehat{q}) 
      ,
\end{align*} 
  which
is correct also for ${\widehat{E} = 0}$ as ${\widehat{x}(0)=0}$. Formula for ${\Newopt\Big(\widehat{x}( \widehat{E})\Big)}$ comes as the~obvious consequence: 
\begin{align*}
\Newopt\Big(\widehat{x}( \widehat{E})\Big) =  
   -\mu_1 + \quantAA \sigma_1 + \widehat{E} \left( \frac{\quantAA \beta_C^{} }{ \alpha_C^{}}- 1\right) + \frac{|\widehat{E}| }{ \alpha_C^{}} \sqrt{ \Delta }.
 \end{align*}
 Now we find~$\CoVaR$-efficient portfolios. The only ones that might satisfy the~required conditions are ${x(\widehat{E}) = (1- \mathbbm{1}_{n-1}^{\mathrm{T}} \widehat{x}(\widehat{E}), \widehat{x}(\widehat{E})^\mathrm{T})^\mathrm{T}} $ portfolios as graph of~${g(\widehat{E}) := \Newopt\Big(\widehat{x}( \widehat{E})\Big)}$ is the~lower boundary of~${\lbrace  \big(x^\mathrm{T} \mu, \CoVaR(x) \big)\mid x \in \mathbb{R}^n \rbrace}$. Function $g(\widehat{E})$ is a~continuous piecewise function comprising two linear functions. It can be easily observed that  whether for a~given $\widehat{E}$ portfolio ${x(\widehat{E})}$ is $\CoVaR$-efficient depends solely on the ratio of ${\quantAA \beta_C^{} \slash \alpha_C^{}- 1}$ and ${\alpha_C^{} \sqrt{ \Delta }}$, or, to be more specific, on~the~inequalities between ${\quantA \beta_C^{} - \alpha_C^{}}$,${-\sqrt{\Delta}}$ and ${\sqrt{\Delta}}$. 

\subsection{Additional remarks}
First note that without assuming linear independence of~${\mu, \mathbbm{1}_n, q}$ there is ${q = \xi_1 \mu + \xi_2 \mathbbm{1}_n}$ (by previous assumption ${\mu}$ and ${\mathbbm{1}_n}$ are linearly independent). Therefore the optimization problem \eqref{CoVaR=3} would have the same critical set as that of Markowitz.

Observe also that for ${\quantA \rightarrow 0}$ function $\widehat{x}(\widehat{E})$ converges to a linear function which is only to be expected by looking at~the~function defined in~\eqref{CoVaR=3}. Still, the problem is not equivalent to that of minimizing $\sigma$ (or $\VaRa$) as  ${\sqrt{x^{\mathrm{T}} Q x}=\sqrt{\widehat{x}^{\mathrm{T}} \widehat{Q} \widehat{x}}}$, not ${\sqrt{x^{\mathrm{T}}\Sigma x}}$, is minimized.
 A linear function is achieved in~no~other way, as ${\widehat{\mu} \nparallel \widehat{q}}$ due to linear independence of~${\mu, \mathbbm{1}_n, q}$.
  Therefore, the image of the \mbox{\lq $\operatorname{CoVaR}^{=}$-critical polyline\rq} \ in~${E \mapsto \Big(E, \Newopt\Big(\widehat{x}( \widehat{E})\Big)\Big)}$ consist of~two rays and with our assumption concerning $\Delta$ is never a line.

Now let ${\Delta=0} $ (\mbox{i.e.} ${s=1}$).  Note here that ${p = \beta_C^{} \slash \alpha_C^{} \widehat{E}}$ from the~lemma is not necessarily a~positive number.
Should we solve:
\begin{equation*}
 \left\lbrace \begin{array}{lcccccc}
			  \widehat{x}^{\mathrm{T}} \widehat{\mu} = \widehat{E}\\
			  \widehat{x}^{\mathrm{T}} \widehat{q} = \widehat{t}\\
\end{array} \right. 
\end{equation*}
for any solution ${\widehat{x}_{\widehat{t}}(\widehat{E})}$ we would get ${\funH \big(\widehat{x}_{\widehat{t}} (\widehat{E})\big) \overset{\widehat{t} \rightarrow -\infty }{\longrightarrow} p}$ (solution being unique for given ${\widehat{t},\widehat{E}}$ only in case of~${n=3}$).

\section{Future research} 
Present work but lightly touches the wide and complex subject~of portfolio optimization for $\CoVaR$. For any question answered 
 few more are raised. What if the normality assumption was to be dropped? What if the Gauss distribution was to be replaced by another one? Will the results hold for~$\operatorname{CoVaR}$?  How solving the problem for various families of copulas, as done in~\cite{bernardi2016covar}, would change the outcome?
Also,~\cite{MaiScha} show that $\CoVaR$ is not monotonic with respect to $\rho_{X,Y}$---how badly does it affect the presented model?

\section*{Acknowledgments}
 
The author thanks Professor Piotr Jaworski for his invaluable advice and insightful comments.

This research did not receive any specific grant from funding agencies in the public, commercial, or \mbox{not-for-profit sectors}.


\bibliography{VaR_CoVaR_bibfile}

\end{document}